\documentclass[11pt]{amsart}
\usepackage{amssymb, latexsym}
\theoremstyle{plain}
\newtheorem{theorem}{Theorem}

\newtheorem*{thm-cheb}{Theorem (Chebyshev)}

\newtheorem{proposition}{Proposition}

\newtheorem*{2'}{Theorem 2'}
\newtheorem*{3'}{Theorem 3'}

\theoremstyle{remark}

\newtheorem*{Remark 1}{Remark 1}
\newtheorem*{Remark 2}{Remark 2}
\newtheorem*{Remark 3}{Remark 3}
\newtheorem*{Remark 4}{Remark 4}

\numberwithin{equation}{section}

\begin{document}

\title[ Inversions   for distribution-biased/shifted  permutations]
 {Comparing the inversion statistic for distribution-biased and distribution-shifted
  permutations with the geometric and the GEM distributions}

\author{Ross G. Pinsky}


\address{Department of Mathematics\\
Technion---Israel Institute of Technology\\
Haifa, 32000\\ Israel}
\email{ pinsky@math.technion.ac.il}

\urladdr{http://www.math.technion.ac.il/~pinsky/}

\date{}

\begin{abstract}
Given a probability distribution $p:=\{p_k\}_{k=1}^\infty$ on the positive integers, there are two natural ways to construct a random permutation in $S_n$ or a random permutation of $\mathbb{N}$ from IID samples from $p$.
One is  called the $p$-biased construction and the other  the $p$-shifted construction.
 In the first part of the paper we consider the case that the distribution $p$ is the
geometric distribution with parameter $1-q\in(0,1)$.  In this case, the $p$-shifted random permutation has the Mallows distribution with parameter $q$.
Let
 $P_n^{b;\text{Geo}(1-q)}$ and $P_n^{s;\text{Geo}(1-q)}$
 denote the biased and the shifted distributions on $S_n$.
The expected number of inversions of a permutation under $P_n^{s;\text{Geo}(1-q)}$ is greater than   under $P_n^{b;\text{Geo}(1-q)}$,
and under either of these distributions, a permutation  tends to have many fewer inversions than it would have
  under the uniform distribution.
 For fixed $n$, both  $P_n^{b;\text{Geo}(1-q)}$ and $P_n^{s;\text{Geo}(1-q)}$ converge weakly as $q\to1$ to the uniform distribution on $S_n$.
 We compare the biased and the shifted distributions by studying the inversion statistic under
 $P_n^{b;\text{Geo}(q_n)}$ and $P_n^{s;\text{Geo}(q_n)}$ for various rates of convergence of $q_n$ to 1.
In the second part of the paper we consider  $p$-biased and $p$-shifted permutations for the case that the distribution $p$ is itself random and distributed as a GEM$(\theta)$-distribution.
In particular,  in both the GEM$(\theta)$-biased and the GEM$(\theta)$-shifted cases,
the expected number of inversions behaves asymptotically as it does
 under the Geo$(1-q)$-shifted distribution with $\theta=\frac q{1-q}$.
This allows one to consider the GEM$(\theta)$-shifted case as the random counterpart of the Geo$(q)$-shifted case.
We also consider another $p$-biased distribution with random $p$    for which the expected number  of inversions  behaves asymptotically
as it does under the Geo$(1-q)$-biased case with $\theta$ and $q$ as above, and with
 $\theta\to\infty$ and $q\to1$.
\end{abstract}

\maketitle
\section{Introduction and Statement of Results}
A permutation of $\mathbb{N}$ is a 1-1  map from $\mathbb{N}$ onto itself.
Let $p:=\{p_k\}_{k=1}^\infty$ be a probability distribution on the positive integers, with $p_k>0$ for all $k$. From this distribution, we describe two methods for creating a random permutation $\Pi:=\{\Pi_k\}_{k=1}^\infty$ of $\mathbb{N}$.
Take a countable sequence of independent samples  from the distribution $p$: $n_1,n_2,\cdots $. The first method is to  define $\Pi_k$ to be the $k$th distinct number to appear in the sequence $\{n_1,n_2,\cdots\}$.
Thus, for example, if the sequence of independent samples from $p$ is $7,3,4,3,7,2,5,\cdots$, then the permutation $\Pi$ begins with $\Pi_1=7, \Pi_2=3,\Pi_3=4,\Pi_4=2,\Pi_5=5$.
Such a random permutation is called a \it $p$-biased permutation.\rm\
The second method is defined as follows. Let $\Pi_1=n_1$ and then for $k\ge2$, let
  $\Pi_k=\psi_k(n_k)$, where $\psi_k$ is the increasing bijection from $\mathbb{N}$ to $\mathbb{N}-\{\Pi_1,\cdots,\Pi_{k-1}\}$.
Thus, the sequence of samples  $7,3,4,3,7,2,5,\cdots$ yields the permutation $\Pi$ beginning with
$\Pi_1=7,\Pi_2=3,\Pi_3=5, \Pi_4=4,\Pi_5=11,\Pi_6=2,\Pi_7=10$. Such a permutation is called a \it $p$-shifted \rm\ permutation.

For any fixed $n\in\mathbb{N}$, one can also obtain a $p$-biased or a  $p$-shifted random permutation of $[n]:=\{1,\cdots, n\}$,
which we denote by $\Pi^{(n)}=\{\Pi^{(n)}_k\}_{k=1}^n$.
 Indeed, we simply
ignore all values that land outside of $[n]$ and stop the process after a finite number of steps, when every value in $[n]$ is obtained. Thus, for example, if we take $n=5$, and if, as before, we sample  the sequence
$7,3,4,3,7,2,5,\cdots$, then we obtain the permutation $34251\in S_5$ in the biased case and $35421\in S_5$ in the shifted case.

Let $P_\infty^{b;\{p_k\}}$ and $P_\infty^{s;\{p_k\}}$
denote the biased and shifted distributions on the permutations of $\mathbb{N}$, induced by the random permutation
$\Pi$,
 and let
$P_n^{b;\{p_k\}}$ and $P_n^{s;\{p_k\}}$ denote the biased and shifted distributions on $S_n$, the set of permutations of
$[n]$, induced by the random permutation $\Pi^{(n)}$.
It is easy to see from the construction   that
$P_n^{b;\{p_k\}}$ and $P_n^{s;\{p_k\}}$ converge weakly to
$P_\infty^{b;\{p_k\}}$ and $P_\infty^{s;\{p_k\}}$ as $n\to\infty$, in the sense that for each $j\in\mathbb{N}$, one has
$$
\begin{aligned}
&P_\infty^{b;\{p_k\}}\big((\sigma_1,\cdots, \sigma_j)\in\cdot\big)=\lim_{n\to\infty}P_n^{b;\{p_k\}}\big((\sigma_1,\cdots, \sigma_j)\in \cdot\big);\\
&P_\infty^{s;\{p_k\}}\big((\sigma_1,\cdots, \sigma_j)\in\cdot\big)=\lim_{n\to\infty}P_n^{s;\{p_k\}}\big((\sigma_1,\cdots, \sigma_j)\in \cdot\big),
\end{aligned}
$$
where $\sigma=\sigma_1\sigma_2\cdots$ denotes a canonical permutation of $\mathbb{N}$, and $\sigma=\sigma_1\cdots,\sigma_n$ denotes a canonical permutation in $S_n$.

In this paper,
we  study the behavior of the inversion statistic.
We first consider $p$-biased and $p$-shifted random permutations in the case that the distribution $p$ is the geometric distribution Geo$(1-q)$:
\begin{equation}\label{p}
p_k=(1-q)q^{k-1},\ k=1,2,\cdots,
\end{equation}
where $q\in(0,1)$.
Then we consider $p$-biased and $p$-shifted random permutations in the case
that the distribution $p$ is itself random and distributed according to the GEM$(\theta)$ distribution, for $\theta>0$.
As will be seen, in the $p$-shifted situation, but not in the $p$-biased situation,  the GEM$(\theta)$ case
may be thought of as a natural random counterpart of the  deterministic Geo$(1-q)$ case, with $q$ and $\theta$ related
by $q=\frac\theta{\theta+1}$ or equivalently, $\theta=\frac q{1-q}$.
This leads us to also consider an alternative  random distribution in the $p$-biased case that can better be considered
as the natural random counterpart of the Geo$(1-q)$ case, with $q$ and $\theta$ related as above.

We begin with the Geo$(1-q)$-biased and Geo$(1-q)$-shifted random permutations.
Denote the corresponding biased and shifted distributions on the permutations of $\mathbb{N}$ and on $S_n$ by
$P_\infty^{b;\text{Geo}(1-q)}, P_\infty^{s;\text{Geo}(1-q)}, P_n^{b;\text{Geo}(1-q)}, P_n^{s;\text{Geo}(1-q)}$.
It is known \cite{GO} that $P_n^{s;\text{Geo}(1-q)}$, the Geo$(1-q)$-shifted distribution on $S_n$,
is actually the \it Mallows \rm distribution with parameter $q$. The Mallows distribution with parameter $q$ is the probability
measure on  $S_n$
that assigns to  each permutation $\sigma\in S_n$ a probability proportional to $q^{\mathcal{I}_n(\sigma)}$, where
$\mathcal{I}_n(\sigma)$ is the number of inversions in $\sigma$; that is
$\mathcal{I}_n(\sigma)=\sum_{1\le i<j\le n}1_\{\sigma_j<\sigma_i\}$.
We extend the inversion statistic $\mathcal{I}_n$ to permutations $\sigma=\sigma_1\sigma_2\cdots$ of $\mathbb{N}$ by defining
$$
\mathcal{I}_n(\sigma)=\sum_{1\le i<j\le n}1_{\{\sigma^{-1}_j<\sigma^{-1}_i\}}=\sum_{\stackrel{1\le k<l<\infty} {\sigma_k,\sigma_l\le n}}1_{\{\sigma_l<\sigma_k\}}.
$$
\noindent \bf Remark.\rm\ From the constructions above, it follows immediately that
the distribution of $\mathcal{I}_n$ under
 $P_\infty^{b;\text{Geo}(1-q)}$ coincides with its distribution  under $P_n^{b;\text{Geo}(1-q)}$,
 and the distribution of $\mathcal{I}_n$
under $P_\infty^{s;\text{Geo}(1-q)}$ coincides with its distribution under under $P_n^{s;\text{Geo}(1-q)}$.
Thus in the sequel, asymptotic results concerning the behavior of  $\mathcal{I}_n$ 
under $P_n^{b;\text{Geo}(1-q)}$ or $P_n^{s;\text{Geo}(1-q)}$ will be stated using the fixed probability
measure $P_\infty^{b;\text{Geo}(1-q)}$ or  $P_\infty^{s;\text{Geo}(1-q)}$.

 We will prove the following proposition.
 \begin{proposition}\label{dom}
 For all $1\le i<j<\infty$,  $1_{\{\sigma^{-1}_j<\sigma^{-1}_i\}}$
under $P_\infty^{s;\text{Geo}(1-q)}$ stochastically dominates $1_{\{\sigma^{-1}_j<\sigma^{-1}_i\}}$    under $P_\infty^{b;\text{Geo}(1-q)}$. The domination is strict if $j-i\ge2$.
 \end{proposition}
From the proposition and the linearity of the expectation it is immediate that
\begin{equation}\label{stochdomexp}
E_\infty^{b;\text{Geo}(1-q)}\mathcal{I}_n>E_\infty^{s;\text{Geo}(1-q)}\mathcal{I}_n,\ \text{for}\ n\ge3.
\end{equation}

It is easy to see from the construction that as $q\in(0,1)$ approaches 1, both the Geo$(1-q)$-biased distribution $P_n^{b;\text{Geo}(1-q)}$ and the Geo$(1-q)$-shifted  distribution $P_n^{s;\text{Geo}(1-q)}$ converge weakly to
the uniform measure on $S_n$.
We compare the behavior of the inversion statistic $\mathcal{I}_n$ under
$P_n^{b;\text{Geo}(1-q)}$ and $P_n^{s;\text{Geo}(1-q)}$
(or equivalently, under $P_\infty^{b;\text{Geo}(1-q)}$ and $P_\infty^{s;\text{Geo}(1-q)}$, by the remark after
Proposition \ref{dom})
for   various rates of convergence of $q_n$ to 1.
We begin however with the case of fixed $q\in(0,1)$.
The notation $w-\lim_{n\to\infty}$ will be used to denote convergence in distribution of a sequence of random variables.
\begin{proposition}\label{fixed}
\noindent  Let $q\in(0,1)$.

\noindent i.
\begin{equation}\label{expfixedq-biased}
\lim_{n\to\infty}\frac{E_\infty^{b;\text{Geo}(1-q)}\mathcal{I}_n}n=\sum_{k=1}^\infty\frac1{1+q^{-k}},
\end{equation}
and
\begin{equation}\label{qto1}
\lim_{q\to1}(1-q)\lim_{n\to\infty}\frac{E_\infty^{b;\text{Geo}(1-q)}\mathcal{I}_n}n=\log 2.
\end{equation}
Furthermore, under $P_\infty^{b;\text{Geo}(1-q)}$,
\noindent $\text{\rm w-lim}_{n\to\infty}\frac{\mathcal{I}_n}n=\sum_{k=1}^\infty\frac1{1+q^{-k}}$.

\noindent ii.
 \begin{equation}\label{expfixedq-shifted}
\lim_{n\to\infty}\frac{E_\infty^{s;\text{Geo}(1-q)}\mathcal{I}_n}n=\frac q{1-q},
\end{equation}
and
$$
\lim_{q\to1}(1-q)\lim_{n\to\infty}\frac{E_\infty^{s;\text{Geo}(1-q)}\mathcal{I}_n}n=1.
$$
Furthermore,
under $P_\infty^{s;\text{Geo}(1-q)}$,
$\text{\rm w-lim}_{n\to\infty}\frac{\mathcal{I}_n}n=\frac q{1-q}$.

\end{proposition}
\begin{theorem}\label{1}
\noindent a. Let $q_n=1-\frac c{n^\alpha},\ \text{with}\ c>0\ \text{and}\ \alpha\in(0,1)$.

i. Under $P_\infty^{b;\text{Geo}(1-q_n)}$,

$$
\text{\rm w-lim}_{n\to\infty}\frac{\mathcal{I}_n}{n^{1+\alpha}}=
\lim_{n\to\infty}\frac{E_\infty^{b;\text{Geo}(1-q_n)}\mathcal{I}_n}{n^{1+\alpha}}=
\frac{\log 2}c.
$$

ii. Under $P_\infty^{s;\text{Geo}(1-q_n)}$,

$$
\text{\rm w-lim}_{n\to\infty}\frac{\mathcal{I}_n}{n^{1+\alpha}}=
\lim_{n\to\infty}\frac{E_\infty^{s;\text{Geo}(1-q_n)}\mathcal{I}_n}{n^{1+\alpha}}=
\frac1c.
$$

\noindent b.  Let $q_n=1-\frac cn,\ \text{with}\ c>0$.

 i. Under $P_\infty^{b;\text{Geo}(1-q_n)}$,

$$
\text{\rm w-lim}_{n\to\infty}\frac{\mathcal{I}_n}{n^2}=
\lim_{n\to\infty}\frac{E_\infty^{b;\text{Geo}(1-q_n)}\mathcal{I}_n}{n^2}=
\frac1{c^2}
\int_0^{1-e^{-c}}\frac{\log(1-\frac x2)}{x-1}dx:=I_b(c).
$$

ii. Under $P_\infty^{s;\text{Geo}(1-q_n)}$,

$$
\text{\rm w-lim}_{n\to\infty}\frac{\mathcal{I}_n}{n^2}=
\lim_{n\to\infty}\frac{E_\infty^{s;\text{Geo}(1-q_n)}\mathcal{I}_n}{n^2}=
\frac1{c^2}
\int_0^{1-e^{-c}}\big(\frac1{1-x}+\frac{\log(1-x)}x\big)dx:=I_s(c).
$$
\noindent Also, $I_b(c)<I_s(c)$,
$\lim_{c\to\infty}I_b(c)=\lim_{c\to\infty}I_s(c)=0$ and\newline
 $\lim_{c\to0}I_b(c)=\lim_{c\to0}I_s(c)=\frac14$.
\

\

\noindent c. Let $q_n=1-o(\frac 1n)$. Then
$$
\lim_{n\to\infty}\frac{E_\infty^{b;\text{Geo}(1-q)}\mathcal{I}_n}{n^2}=
\lim_{n\to\infty}\frac{E_\infty^{s;\text{Geo}(1-q)}\mathcal{I}_n}{n^2}=
\frac14.
$$
\end{theorem}
\bf\noindent Remark.\rm\  The  dominance in expectation of the inversion statistic under
$P_n^{s;\text{Geo}(1-q_n)}$ as compared to under  $P_n^{b;\text{Geo}(1-q_n)}$  disappears asymptotically
if $q_n=1-o(\frac1n)$. Indeed, in such a case, both distributions mimic the uniform distribution for which
it is well-known that $\lim_{n\to\infty}\frac{E\mathcal{I}_n}{n^2}=\frac14$.
\medskip

We now consider $p$-biased and $p$-shifted random permutations in the case that the distribution $p$ is itself random and distributed according
to the GEM$(\theta)$ distribution, which we now describe. Let $\{W_k\}_{k=1}^\infty$ be IID random variables taking values
in $(0,1)$. Define a random sequence $\{\mathcal{P}_k\}_{k=1}^\infty$, deterministically satisfying $\sum_{k=1}^\infty \mathcal{P}_k=1$,  by
\begin{equation}\label{RAM}
\mathcal{P}_1=W_1,\ \  \mathcal{P}_k=(1-W_1)\cdots (1-W_{k-1})W_k, \ \  k\ge 2.
\end{equation}
 Such a random distribution is called a
\it random allocation model (RAM)\rm\ or a \it stick-breaking model.\rm\ The GEM$(\theta)$ distribution with $\theta>0$ is the RAM model
in the case that the IID sequence    $\{W_k\}_{k=1}^\infty$ has the Beta$(1,\theta)$-distribution; namely the distribution
with density $\theta(1-w)^{\theta-1},\ 0<w<1$.

We denote by $P_\infty^{b;\text{\rm GEM}(\theta)}$ and  $P_\infty^{s;\text{\rm GEM}(\theta)}$  respectively the corresponding biased and shifted  distributions on permutations of $\mathbb{N}$,
and call them  the GEM$(\theta)$-biased  and the GEM$(\theta)$-shifted distributions.
 Note
that we are in the annealed setting. That is, we  sample a sequence $\{p_k\}_{k=1}^\infty$ from the GEM$(\theta)$-distributed
random variables $\{\mathcal{P}_k\}_{k=1}^\infty$ and use this realization to construct a $p$-biased and a $p$-shifted random permutation
of $\mathbb{N}$.  We have
$$
P_\infty^{*;\text{\rm GEM}(\theta)}(\thinspace\cdot\thinspace)=\int P_\infty^{*;\{p_k\}}(\thinspace\cdot\thinspace)dP_{\theta}\{\mathcal{P}_k\}=\{p_k\}),\ \text{for}\ *=b\ \text{or}\ *=s,
$$
 where $P_{\theta}$ is the
probability measure on  the GEM$(\theta)$-distributed sequence $\{\mathcal{P}_k\}_{k=1}^\infty$.
(With an abuse of notation, we will also use $P_\theta$ to denote the probability measure associated with the sequence $\{W_k\}_{k=1}^\infty$ of IID Beta$(1,\theta)$-distributed random variables used to construct the sequence
$\{\mathcal{P}_k\}_{k=1}^\infty$.)
In the same way as in the deterministic case, we can also define
 $P_n^{b;\text{\rm GEM}(\theta)}$ and  $P_n^{s;\text{\rm GEM}(\theta)}$
 on $S_n$. Analogous to the deterministic case, $\mathcal{I}_n$ has the same distribution under
 $P_n^{b;\text{\rm GEM}(\theta)}$ or $P_n^{s;\text{\rm GEM}(\theta)}$
as it does under $P_\infty^{b;\text{\rm GEM}(\theta)}$ or  $P_\infty^{s;\text{\rm GEM}(\theta)}$.

For the Beta$(1,\theta)$-distributed IID random variables $\{W_k\}_{k=1}^\infty$,
we have $E_\theta W_1=\frac1{\theta+1}$ and therefore $E_\theta(1-W_1)=\frac\theta{1+\theta}$.
Thus, comparing the random distribution on $\mathbb{N}$ given by a realization of $\{\mathcal{P}_k\}_{k=1}^\infty$ as in \eqref{RAM}, with $\{W_k\}_{k=1}^\infty$ as above,
with the deterministic geometric distribution on $\mathbb{N}$ given in \eqref{p}, it is  natural to compare
the  Geo$(1-q)$-biased or shifted distribution  to the GEM$(\theta)$-biased or shifted distribution, with
$q$ and $\theta$ related by $q=\frac\theta{\theta+1}$, or equivalently, $\theta=\frac q{1-q}$.
It turns out that with respect to the inversion statistic, this comparison is apt in the shifted case, but not in the biased case. We will prove the following results.
\begin{theorem}\label{extra}
Let $\theta>0$.  For  $P_\theta$-almost all $\{\mathcal{P}_k\}_{k=1}^\infty=\{p_k\}_{k=1}^\infty$,
\begin{equation}\label{weakgemshift}
w-\lim_{n\to\infty}\frac{\mathcal{I}_n}n=\sum_{k=1}^\infty k\thinspace \mathcal{P}_{k+1}=\sum_{k=1}^\infty k\thinspace W_{k+1}\prod_{i=1}^k(1-W_i),
\end{equation}
where $w-\lim_{n\to\infty}$ denotes the weak limit under the measure $P_\infty^{s;\{p_k\}}$.
Furthermore,
\begin{equation}\label{expextra}
\lim_{n\to\infty}\frac{E_\infty^{s;\text{\rm GEM}(\theta)}\mathcal{I}_n}n=\theta.
\end{equation}
\end{theorem}
\begin{theorem}\label{biasedgemthm}
Let $\theta>0$.
Then
\begin{equation}\label{expbiasedgem}
\lim_{n\to\infty}\frac{E_\infty^{b;\text{\rm GEM}(\theta)}\mathcal{I}_n}n=\theta.
\end{equation}
\end{theorem}
\bf\noindent Remark 1.\rm\  The calculations involved in the proof of Theorem \ref{biasedgemthm} are the most interesting ones in the paper, and contain several twists and novelties.

\medskip
\noindent \bf Remark 2.\rm\
In light of \eqref{stochdomexp},
 it is not surprising   that the right hand side of \eqref{expfixedq-shifted}
is larger than the right hand side of \eqref{expfixedq-biased}. Note however that
the right hand sides of \eqref{expextra} and \eqref{expbiasedgem}  are the same.

\medskip

 With regard to the discussion in the paragraph preceding Theorem \ref{extra},
 compare \eqref{expextra} to \eqref{expfixedq-shifted}. From this,     in the shifted case
 $P_\infty^{s;\text{\rm GEM}(\theta)}$ might be thought of as the natural random counterpart of
 $P_\infty^{s;\text{Geo}(1-q)}$, 
 with $\theta=\frac q{1-q}$.
 However, comparing  \eqref{expbiasedgem} to \eqref{expfixedq-biased} shows that such a connection does not carry over to
$P_\infty^{b;\text{\rm GEM}(\theta)}$ and $P_\infty^{b;\text{Geo}(1-q)}$
in the biased case.
In light of this, we now consider another family of $p$-biased  distributions with random distribution $p$
which, as we shall see,  better deserves to be considered as the natural
random counterpart to the family of  $P^{b;\text{Geo}(1-q)}$-distributions.
Let $\{U_k\}_{k=1}^\infty$ be a sequence of IID random variables  distributed uniformly on $[0,1]$.
Denote expectation with respect to these random variables by the generic $E$.
Let $\theta>0$. Define a random sequence $\{\mathcal{P}_k'\}_{k=1}^\infty$
by
$$
\mathcal{P}_k'=\prod_{i=1}^k U_i^\frac1\theta.
$$
Let
$$
D=\sum_{k=1}^\infty \mathcal{P}_k'=\sum_{k=1}^\infty \prod_{i=1}^kU_i^\frac1\theta,
$$
and define the random sequence $\{\mathcal{P}_k\}_{k=1}^\infty$,  deterministically satisfying
$\sum_{k=1}^\infty \mathcal{P}_k=1$,  by
\begin{equation}\label{Pindep}
\mathcal{P}_k=\frac{\mathcal{P}_k'}{D}=\frac1{D}\prod_{i=1}^k U_i^\frac1\theta.
\end{equation}
We consider the $p$-biased distribution with $p$ distributed as $\{\mathcal{P}_k\}_{k=1}^\infty$, and denote this
distribution by $P_\infty^{b;\text{\rm IID-prod}(\theta)}$.
We note that the normalization random variable $D$ is known to have the so-called generalized Dickman
distribution with parameter $\theta$ \cite{P18}. However, from the construction,  $D$ does not enter into the formulas for the inversion probabilities; for example,
\begin{equation*}\label{Pprimeinver}
P_\infty^{b;\text{\rm IID-prod}(\theta)}(\sigma^{-1}_j<\sigma^{-1}_i)=E\frac{\mathcal{P}_j'}{\mathcal{P}_i'+\mathcal{P}_j'}.
\end{equation*}
Note that $U_k^\frac1\theta$ has density $\theta x^{\theta-1}, x\in[0,1]$; thus $U_k^\frac1\theta\stackrel{\text{dist}}{=}1-W_k$, where $W_k$ has the
Beta$(1,\theta)$ distribution. In particular, $EU_k^\frac1\theta=\frac\theta{\theta+1}$.
Thus, letting $\{W_k\}_{k=1}^\infty$ be an IID sequence of Beta$(1,\theta)$-distributed random variables,
the random sequence $\{\mathcal{P}_k\}_{k=1}^\infty$ constructed above in \eqref{Pindep} can also be constructed in the following
 equivalent manner:
$$
\mathcal{P}_k'=\prod_{i=1}^k (1-W_i);
$$
$$
D=\sum_{k=1}^\infty \mathcal{P}_k'=\sum_{k=1}^\infty \prod_{i=1}^k(1-W_i);
$$
\begin{equation}\label{Pindepagain}
\mathcal{P}_k=\frac{\mathcal{P}_k'}{D}=\frac1{D}\prod_{i=1}^k(1-W_i).
\end{equation}
Comparing \eqref{p}, \eqref{RAM} and \eqref{Pindep} (or \eqref{Pindepagain}), we  suggest that, with $\theta$ and $q$ related by $q=\frac\theta{\theta+1}$, or equivalently, $\theta=\frac{q}{1-q}$, the distribution
$P_\infty^{b;\text{\rm IID-prod}(\theta)}$ rather than the distribution
$P_\infty^{b;\text{\rm GEM}(\theta)}$ should be considered as the natural random counterpart of
the distribution $P_\infty^{b;\text{\rm Geo}(q)}$, at least as $q\to1$ and $\theta\to\infty$.
The following theorem supports this claim; indeed,
compare \eqref{qto1} to
\eqref{expbiasedgem} and \eqref{expbiasedindep}.
\begin{theorem}\label{indeptheta}
Let $\theta>0$.
Then
\begin{equation}\label{expbiasedindep}
\lim_{n\to\infty}\frac{E_\infty^{b;\text{\rm IID-prod}(\theta)}\mathcal{I}_n}n=\theta\log 2.
\end{equation}
\end{theorem}

\medskip

Note that for the shifted case in Theorem \ref{extra} we have a weak law of large numbers as well as an asymptotic result for the expected value, whereas for the biased case in Theorems \ref{biasedgemthm} and \ref{indeptheta} we only have
an asymptotic result for the expected value. The following proposition, of independent interest, concerning the generic shifted case constructed from an  arbitrary deterministic
distribution on $\mathbb{N}$,
 makes it easier to prove a weak law in the shifted case.
The proposition will also be used in the proof of the law of large numbers for the shifted case in Proposition \ref{fixed} and Theorem \ref{1}.
\begin{proposition}\label{indepinver}
Let $p:=\{p_k\}_{k=1}^\infty$ be a probability distribution on $\mathbb{N}$, and let $P_\infty^{s;\{p_k\}}$ denote the corresponding $p$-shifted distribution on the permutations of $\mathbb{N}$.
Let $I_{<j}(\sigma)$ denote the number of inversions involving  pairs $\{\{i,j\}:1\le i<j\}$, for
 $\sigma$ a permutation of $\mathbb{N}$.
Under $P_\infty^{s;\{p_k\}}$, the random variables
$\{I_{<j}\}_{j=1}^\infty$ are independent. Furthermore, the distribution of $I_{<j}$ is given by
\begin{equation}\label{<jdist}
P_\infty^{s;\{p_k\}}(I_{<j}=l)=\frac{p_{l+1}}{\sum_{k=1}^jp_k},\ l=0,1,\cdots, j-1.
\end{equation}
\end{proposition}
\noindent \bf Remark 1.\rm\ In the case that the distribution $p$ is the Geo$(1-q)$ distribution, the proposition shows
that $I_{<j}$ is distributed as a truncated geometric distribution with parameter
$1-q$,
starting from 0 and truncated at $j-1$: $P_\infty^{s;\text{Geo}(1-q)}(I_{<j}=l)=\frac{(1-q)q^l}{1-q^j}$,\ $l=0,1,\cdots, j-1$.
Actually, Proposition \ref{indepinver} in the case that $p$ is  the  Geo$(1-q)$ distribution is well-known and follows from an alternative construction of the Mallows distribution--see  \cite{P} for example. This alternative construction appears generically
in Remark 2 below.

\noindent \bf Remark 2.\rm\ From Proposition \ref{indepinver} it follows that the $p$-shifted random permutation
$\Pi^{(n)}$ (or $\Pi$)  can be constructed in an alternative manner by sequentially placing the numbers
$\{1,\cdots, n\}$ (or $\{1,2,\cdots\}$)  down on a line
at various positions between the numbers that have already been placed down. First place down the number 1.
For $j\ge2$, assume that the  numbers $\{1,\cdots, j-1\}$ have already been placed down.
Then there are  $j$ possible spaces in which to place the number $j$; namely, to the right of any of the $j-1$ numbers that have already been placed down, or to the left of the leftmost number that has already been placed down.
For $l=0,\cdots, j-1$,
with probability $\frac{p_{l+1}}{\sum_{k=1}^jp_k}$
place the number $j$ in the $(l+1)$-th rightmost position. Note that this gives
$1_{<j}=l$.

\medskip

Although we won't need it here, we note that four out of the five models of
random permutations discussed above are examples of  \it strictly regenerative permutations.\rm\
(The exception is the GEM($\theta$)-shifted case.)
For a permutation $\pi=\pi_{a+1}\pi_{a+2}\cdots\pi_{a+m}$, of $\{a+1,a+2,\cdots, a+m\}$, define
 $\text{red}(\pi)$, the reduced permutation of $\pi$, to be the permutation in $S_m$ given by
 $\text{red}(\pi)_i=\pi_{a+i}-m$.
 A random permutation  is strictly regenerative if
 for almost every realization $\Pi$ of the random permutation, there exist
$0=T_0< T_1<T_2<\cdots$ such that $\Pi([T_j])=[T_j],\ j\ge1$,  and $\Pi([m])\neq[m]$ if $m\not\in\{T_1,T_2,\cdots\}$,  and such that the
random variables $\{T_k-T_{k-1}\}_{k=1}^\infty$ are IID and the  random permutations
$\{\text{red}(\Pi|_{[T_k]-[T_{k-1}]}\}_{k=1}^\infty$ are IID.
The intervals $\{T_k-T_{k-1}\}_{k=1}^\infty$ are called the \it blocks\rm\ of the  permutation.
The four aforementioned  models
are  \it positive recurrent\rm, which means that
the block length has finite expected value; that is, $ET_1<\infty$. For more on this, see \cite{PT} and references therein.
In particular, in the specific context of Mallows distributions, for fixed $q$, see \cite{GO} for more on general constructions, and see  \cite{BB} for an analysis of the length of the longest increasing subsequence;
for $q_n\to1$, see \cite{BP} for an analysis of the length of the longest increasing subsequence and see  \cite{GP} for an analysis of
the cycle structure.

In section \ref{2-sec} we prove Propositions \ref{dom} and \ref{indepinver}. In section \ref{expectation-sec} we analyze
the expected number of inversions, $E_\infty^{b;\text{Geo}(1-q_n)}I_n$ and   $E_\infty^{s;\text{Geo}(1-q)}I_n$, for $q_n\equiv q$ as in
Proposition \ref{fixed} and for the various cases of $q_n$ as in  Theorem \ref{1}.
In section \ref{wlln-sec},   applications of the second moment method along with the results of section \ref{expectation-sec} yield the proofs of
Proposition \ref{fixed} and Theorem \ref{1}. The proof of Theorem \ref{extra} is given in section \ref{extrastuff},
the proof of Theorem \ref{biasedgemthm} is given in section \ref{GEMofaproof} and the proof of Theorem \ref{indeptheta}
is given in section \ref{proofindeptheta}.

\section{Proofs of Propositions \ref{dom} and \ref{indepinver}}\label{2-sec}
\it\noindent Proof of Proposition \ref{dom}.\rm\
From the construction of the biased case, it is clear that $P_\infty^{b;\text{Geo}(1-q)}(\sigma^{-1}_j<\sigma^{-1}_i)=\frac{p_j}{p_j+p_i}$.
This probability is equal to $\frac{q^j}{q^j+q^i}$.
We now show that $P_\infty^{s;\text{Geo}(1-q)}(\sigma^{-1}_j<\sigma^{-1}_i)\ge\frac{q^j}{q^j+q^i}$.
From the construction of  the shifted case, it is clear that on the first step of the construction, the probability that
$j$ will appear, conditioned on either $i$ or $j$ appearing on that step, is equal to $\frac{p_j}{p_j+p_i}$, which is equal to $\frac{q^j}{q^j+q^i}$.
If the number appearing on the first step is $k\neq i,j$, then the probability that
$j$ will appear on the second step, conditioned on either $i$ or $j$ appearing on that step, depends on the value
of $k$. If $k>j$, then this probability is again $\frac{p_j}{p_j+p_i}=\frac{q^j}{q^j+q^i}$. If $k<i$, then
this probability is $\frac{p_{j-1}}{p_{j-1}+p_{i-1}}=\frac{q^j}{q^j+q^i}$.
However, if $i<k<j$, then this probability is equal to
$\frac{p_{j-1}}{p_{j-1}+p_i}=\frac{q^{j-1}}{q^{j-1}+q^i}> \frac{q^j}{q^j+q^i}$. Thus, the probability
that $j$ will appear on the second step, conditioned on either $i$ or $j$ appearing on that step, and conditioned
on neither of them having already appeared on the first step, is greater or equal
to $\frac{q^j}{q^j+q^i}$, and in fact, strictly greater if $j-i>1$. Continuing in this vein proves the proposition.
\hfill $\square$
\medskip

\noindent \it Proof of Proposition \ref{indepinver}.\rm\
We first prove that the distribution of $1_{<j}$ is given by \eqref{<jdist}.
From the construction of the shifted permutation, it follows that for $i\in\{1,\cdots, j\}$, the
probability that from among the numbers $\{1,\cdots, j\}$, the first one to be placed down
 in the permutation will be $i$ is $\frac{p_i}{\sum_{k=1}^jp_k}$. Thus, in particular, in the case $i=j$, we obtain
 $P_\infty^{s;\{p_k\}}(1_{<j}=j-1)=\frac{p_j}{\sum_{k=1}^jp_k}$.
With probability $\frac{\sum_{k=1}^{j-1}p_k}{\sum_{k=1}^jp_k}$, the number $j$ will not be the first number to be
placed down from among
the numbers $\{1,\cdots, j\}$. It follows from the shifted construction that
conditioned on this event, the probability that
the number $j$ will be the second number to be placed down from among the numbers $\{1,\cdots, j\}$ is equal
to $\frac{p_{j-1}}{\sum_{k=1}^{j-1}p_k}$. Thus, it follows that
$P_\infty^{s;\{p_k\}}(1_{<j}=j-2)=
\frac{\sum_{k=1}^{j-1}p_k}{\sum_{k=1}^jp_k}
\times\frac{p_{j-1}}{\sum_{k=1}^{j-1}p_k}=\frac{p_{j-1}}{\sum_{k=1}^jp_k}$.
Continuing in this vein, we obtain \eqref{<jdist}.

We now prove the independence of the random variables $\{1_{<j}\}_{j=1}^\infty$.
By induction and by what we have already proved, it suffices to show that
\begin{equation}\label{split}
\begin{aligned}
&P_\infty^{s;\{p_k\}}(I_{<2}=a_2,I_{<3}=a_2,\cdots, I_{<j+1}=a_{j+1})=\\
&\frac{p_{a_{j+1}+1}}{\sum_{k=1}^{j+1}p_k}P_\infty^{s;\{p_k\}}(I_{<2}=a_2,I_{<3}=a_2,\cdots, I_{<j}=a_j), \\
& \text{for}\
0\le a_i\le i-1,\ i=2,\cdots, j+1,\ \text{and}\ j\ge2.
\end{aligned}
\end{equation}
As is well known, specifying the values $I_{<2}=a_2,I_{<3}=a_2,\cdots, I_{<j+1}=a_{j+1}$, uniquely
determines a permutation of $\{1,\cdots,  j+1\}$, call it $\sigma=\sigma_1\cdots\sigma_{j+1}$,  specifying
the values $I_{<2}=a_2,I_{<3}=a_2,\cdots, I_{<j}=a_j$, uniquely
determines a permutation of $\{1,\cdots,  j\}$, call it $\tau=\tau_1\cdots \tau_j$, and the permutation obtained
by deleting the  number $j+1$ from $\sigma$ is $\tau$.
Let $i^*=\sigma^{-1}_{j+1}$.
Note then that  $1_{<j+1}(\sigma)=j+1-i^*$. Since we are assuming that
$1_{<j+1}(\sigma)=a_{j+1}$, it follows that $i^*=j+1-a_{j+1}$.

From the observations in the previous paragraph, it  follows from the shifted construction that
\begin{equation}\label{j+1terms}
P_\infty^{s;\{p_k\}}(I_{<2}=a_2,I_{<3}=a_2,\cdots, I_{<j+1}=a_{j+1})=\prod_{i=1}^{j+1}\frac{p_{b_i}}{\sum_{k=1}^{j+2-i}p_k},
\end{equation}
for a certain appropriate choice of $\{b_i\}_{i=1}^{j+1}$, with $1\le b_i\le j+2-i$,
and in particular, $b_{i^*}=j+2-i^*$, and that
\begin{equation}\label{jterms}
P_\infty^{s;\{p_k\}}(I_{<2}=a_2,I_{<3}=a_2,\cdots, I_{<j}=a_j)=\prod_{i=1}^{i^*-1}\frac{p_{b_i}}{\sum_{k=1}^{j+1-i}p_k}
\thinspace\prod_{i=i^*+1}^{j+1}\frac{p_{b_i}}{\sum_{k=1}^{j+2-i}p_k}.
\end{equation}
The difference between the right hand side of
\eqref{j+1terms} and the right hand side of \eqref{jterms} is that
the right hand side of \eqref{j+1terms} has the extra factor $p_{b_{i^*}}$ in its numerator and the extra factor
$\sum_{k=1}^{j+1}p_k$ in its denominator. Now
$\frac{p_{b_{i^*}}}{\sum_{k=1}^{j+1}p_k}=\frac{p_{j+2-i^*}}{\sum_{k=1}^{j+1}p_k}=
\frac{p_{a_{j+1}+1}}{\sum_{k=1}^{j+1}p_k}$.
From these facts,  \eqref{split} follows.
\hfill $\square$

\section{Analysis of the expected number of inversions}\label{expectation-sec}
To calculate the expected number of inversions in the biased case, we write $\mathcal{I}_n=\sum_{1\le i<j\le n}1_{\{\sigma^{-1}_j<\sigma^{-1}_i\}}$. As noted in the proof of Proposition \ref{dom},  $E_\infty^{b;\text{Geo}(1-q)}1_{\{\sigma^{-1}_j<\sigma^{-1}_i\}}=\frac{q^j}{q^j+q^i}$.
Thus
\begin{equation}\label{exp-b}
E_\infty^{b;\text{Geo}(1-q)}\mathcal{I}_n=\sum_{1\le i<j\le n}\frac{q^j}{q^j+q^i}=\sum_{1\le i<j\le n}\frac1{1+q^{i-j}}=
\sum_{k=1}^{n-1}\frac{n-k}{1+q^{-k}}.
\end{equation}

To calculate the the expected number of inversions in the shifted case, we represent $\mathcal{I}_n$ as
$\sum_{j=1}^nI_{<j}$, where $I_{<j}$ is as in Proposition \ref{indepinver}.  By that proposition and Remark 1 following it, we have
$$
\begin{aligned}
&E_\infty^{s;\text{Geo}(1-q)}I_{<j}=\sum_{k=0}^{j-1}\frac{1-q}{1-q^j}kq^k=\frac{(1-q)q}{1-q^j}\sum_{k=0}^{j-1}kq^{k-1}=
\frac{(1-q)q}{1-q^j}\frac{d}{dq}\big(\frac{1-q^j}{1-q}\big)=\\
&\frac {q\big(1+(j-1)q^j-jq^{j-1}\big)}{(1-q^j)(1-q)}.
\end{aligned}
$$
Thus,
$$
E_\infty^{s;\text{Geo}(1-q)}\mathcal{I}_n=\sum_{j=1}^{n-1}\frac {q\big(1+(j-1)q^j-jq^{j-1}\big)}{(1-q^j)(1-q)}.
$$
Performing some algebra \cite{R}, this reduces to
\begin{equation}\label{exp-s}
E_\infty^{s;\text{Geo}(1-q)}\mathcal{I}_n=\frac q{1-q}(n-1)-\sum_{j=1}^{n-1}\frac{jq^j}{1-q^j}.
\end{equation}

We now use \eqref{exp-b} and \eqref{exp-s} to analyze the asymptotic behavior of the expectation
for various choices of $q=q_n$.
\medskip

\noindent \bf The case of fixed $q\in(0,1)$:\rm\

From \eqref{exp-b}, we obtain
\begin{equation}\label{asympexp-b}
\lim_{n\to\infty}\frac{E_\infty^{b;\text{Geo}(1-q)}\mathcal{I}_n}n=\sum_{k=1}^\infty\frac1{1+q^{-k}}.
\end{equation}
Approximating by Riemann sums gives
\begin{equation}\label{approxq}
\int_1^\infty\frac1{1+e^{ax}}dx\le\sum_{k=1}^\infty\frac1{1+q^{-k}}\le\frac q{q+1}+\int_1^\infty\frac1{1+e^{ax}}dx,\ \ a=-\log q.
\end{equation}
We have
\begin{equation}\label{calint}
\int_1^\infty\frac1{1+e^{ax}}dx=\int_1^\infty\frac{e^{-ax}}{e^{-ax}+1}dx=\frac{\log(1+e^{-a})}a=
\frac{\log(1+q)}{-\log q}.
\end{equation}
From \eqref{asympexp-b}-\eqref{calint} it follows
that
\begin{equation}\label{q-b}
\lim_{q\to1}(1-q)\lim_{n\to\infty}\frac{E_\infty^{b;\text{Geo}(1-q)}\mathcal{I}_n}n=\log2.
\end{equation}

From
 \eqref{exp-s} we obtain
 \begin{equation}\label{asympexp-s}
\lim_{n\to\infty}\frac{E_\infty^{s;\text{Geo}(1-q)}\mathcal{I}_n}n=\frac q{1-q}.
\end{equation}

\medskip

\noindent \bf The case of  $q=1-\frac c{n^\alpha},\ \ c>0,\alpha\in(0,1)$.\rm\

From \eqref{exp-b},
we write
\begin{equation}\label{representexp}
E_\infty^{b;\text{Geo}(1-q_n)}\mathcal{I}_n=
n\sum_{k=1}^{n-1}\frac1{1+q_n^{-k}}-\sum_{k=1}^{n-1}\frac k{1+q_n^{-k}}.
\end{equation}
Similar to \eqref{approxq}, we  have
\begin{equation}\label{approxqalpha}
\int_1^n\frac1{1+e^{a_nx}}dx\le \sum_{k=1}^{n-1}\frac1{1+q_n^{-k}}\le\frac {q_n}{q_n+1}+\int_1^{n-1}\frac1{1+e^{a_nx}}dx,\ a_n=-\log q_n.
\end{equation}
Integrating, similar to \eqref{calint}, we obtain
\begin{equation}\label{calintn}
\int_1^n\frac1{1+e^{a_nx}}dx=-\frac1{a_n}\log(1+e^{-a_nx})|_1^n=\frac1{-\log q_n}\big(\log(1+q_n)-\log(1+q_n^n)\big).
\end{equation}
Since $\alpha\in(0,1)$, we have $\lim_{n\to\infty}q_n^n=0$. Thus, from \eqref{approxqalpha} and \eqref{calintn}, the first term on the right hand side of \eqref{representexp} satisfies
\begin{equation}\label{firsttermalpha}
n\sum_{k=1}^{n-1}\frac1{1+q_n^{-k}}\sim\frac{\log 2}c n^{1+\alpha}.
\end{equation}

We now consider the second term on the right hand side of \eqref{representexp}.
 We break it up into two parts. Let $\beta\in(\alpha, \frac{1+\alpha}2)$.
We have
\begin{equation}\label{secondtermalpha}
\sum_{k=1}^{[n^\beta]}\frac k{1+q_n^{-k}}\le n^{2\beta}.
\end{equation}
And we have
\begin{equation}\label{secondagain}
\sum_{[n^\beta]+1}^{n-1}\frac k{1+q_n^{-k}}\le n\sum_{[n^\beta]+1}^{n-1}\frac 1{1+q_n^{-k}}.
\end{equation}
Similar to the argument in \eqref{approxqalpha}-\eqref{firsttermalpha}, we have
\begin{equation}\label{secondyetagain}
\sum_{[n^\beta]+1}^{n-1}\frac 1{1+q_n^{-k}}\sim
\frac1{-\log q_n}\big(\log(1+q_n^{n^\beta})-\log(1+q_n^n)\big)=O(n^\alpha e^{-cn^{\beta-\alpha}}).
\end{equation}
From \eqref{representexp} and \eqref{firsttermalpha}-\eqref{secondyetagain}, we conclude that
\begin{equation}\label{finalnalpha}
\lim_{n\to\infty}\frac{E_\infty^{b;\text{Geo}(1-q_n)}\mathcal{I}_n}{n^{1+\alpha}}=\frac{\log 2}c,\
\ q_n=1-\frac c{n^\alpha},\ \alpha\in(0,1),\  c>0.
\end{equation}

Now we turn to $E_\infty^{s;\text{Geo}(1-q_n)}\mathcal{I}_n$.
From \eqref{exp-s}, we write
\begin{equation}\label{representexps}
E_\infty^{s;\text{Geo}(1-q_n)}\mathcal{I}_n=\frac {q_n}{1-q_n}(n-1)-\sum_{j=1}^{n-1}\frac{jq_n^j}{1-q_n^j}.
\end{equation}
Of course,
\begin{equation}\label{firstterms}
\frac {q_n}{1-q_n}(n-1)\sim \frac{n^{1+\alpha}}c.
\end{equation}
One can check that the function $\frac{xe^{-ax}}{1-e^{-ax}}$ is decreasing for $x\in[1,\infty)$, for $a>0$.
Thus by Riemann sum approximation,
\begin{equation}\label{asymptotically}
\sum_{j=1}^{n-1}\frac{jq_n^j}{1-q_n^j}\sim\int_1^n\frac{xe^{-a_nx}}{1-e^{-a_nx}}dx, \ a_n=-\log q_n.
\end{equation}
We have
\begin{equation}\label{evalint}
\int_1^n\frac{xe^{-a_nx}}{1-e^{-a_nx}}dx=\frac1{a_n^2}\int_{a_n}^{na_n}\frac{ye^{-y}}{1-e^{-y}}dy=
\frac1{(\log q_n)^2}\int_{-\log q_n}^{-n\log q_n}\frac{ye^{-y}}{1-e^{-y}}dy.
\end{equation}
Since $\alpha\in(0,1)$, we conclude from \eqref{asymptotically} and \eqref{evalint} that
\begin{equation}\label{secondterms}
\sum_{j=1}^{n-1}\frac{jq_n^j}{1-q_n^j}\sim\frac{n^{2\alpha}}{c^2}\int_0^\infty\frac{ye^{-y}}{1-e^{-y}}dy.
\end{equation}
From \eqref{representexps}, \eqref{firstterms} and \eqref{secondterms}, we conclude that
\begin{equation}\label{finalalphas}
\lim_{n\to\infty}\frac{E_\infty^{s;\text{Geo}(1-q_n)}\mathcal{I}_n}{n^{1+\alpha}}=\frac1c,\  q_n=1-\frac c{n^\alpha},\ \alpha\in(0,1),\ c>0.
\end{equation}
\medskip

\noindent \bf The case of  $q=1-\frac cn,\ \ c>0$.\rm\

The expectation $E_\infty^{b;\text{Geo}(1-q_n)}\mathcal{I}_n$ is given in \eqref{representexp}.
By Riemann sum approximation,
\begin{equation}\label{alpha1firstb}
\sum_{k=1}^{n-1}\frac{n-k}{1+q_n^{-k}}\sim\int_1^n\frac{n-x}{1+e^{a_nx}}dx,\ \ a_n=-\log q_n.
\end{equation}
Substituting $q_n=1-\frac cn$ in \eqref{calintn}, we obtain
\begin{equation}\label{firstb1asym}
n\sum_{k=1}^{n-1}\frac1{1+q_n^{-k}}\sim \frac{n^2}c\log\frac2{1+e^{-c}}.
\end{equation}

Integrating by parts, we have
\begin{equation}\label{secondterm}
\begin{aligned}
&\int_1^n\frac x{1+e^{a_nx}}dx=\int_1^n\frac{xe^{-a_nx}}{1+e^{-a_nx}}dx=\\
&-\frac x{a_n}\log(1+e^{-a_nx})|_1^n+\frac1{a_n}\int_1^n\log(1+e^{-a_nx})dx.
\end{aligned}
\end{equation}
We have
\begin{equation}\label{firstpart}
\begin{aligned}
&-\frac x{a_n}\log(1+e^{-a_nx})|_1^n=\frac1{-\log q_n}\log(1+q_n)-\frac n{-\log q_n}\log(1+q_n^n)\sim\\
&\frac nc\log2-\frac{n^2}c\log(1+e^{-c})\sim-\frac{n^2}c\log(1+e^{-c}).
\end{aligned}
\end{equation}
Making a change of variables, we have
\begin{equation}\label{secondpart}
\begin{aligned}
&\frac1{a_n}\int_1^n\log(1+e^{-a_nx})dx=\frac1{a_n^2}\int_{e^{-na_n}}^{e^{-a_n}}\frac{\log(1+y)}ydy=\\
&\frac1{(\log q_n)^2}\int_{q_n^n}^{q_n}\frac{\log(1+y)}ydy\sim
\frac{n^2}{c^2}\int_{e^{-c}}^1\frac{\log(1+y)}ydy.
\end{aligned}
\end{equation}
From \eqref{secondterm}-\eqref{secondpart}, we have
\begin{equation}\label{secondb1asym}
\int_1^n\frac x{1+e^{a_nx}}dx\sim n^2\Big(\frac1{c^2}\int_{e^{-c}}^1\frac{\log(1+y)}ydy-\frac1c\log(1+e^{-c})\Big).
\end{equation}
From \eqref{representexp}, \eqref{firstb1asym} and \eqref{secondb1asym}, we
conclude that
\begin{equation}\label{finalnalpha1s}
\begin{aligned}
&\lim_{n\to\infty}\frac{E_\infty^{b;\text{Geo}(1-q_n)}\mathcal{I}_n}{n^2}=\frac1c\log\frac2{1+e^{-c}}+\frac1c\log(1+e^{-c})-
\frac1{c^2}\int_{e^{-c}}^1\frac{\log(1+y)}ydy=\\
&\frac1c\log 2-\frac1{c^2}\int_{e^{-c}}^1\frac{\log(1+y)}ydy=\frac1{c^2}\int_{e^{-c}}^1\big(\frac{\log 2}y-
\frac{\log(1+y)}y\big)dy=\\
&\frac1{c^2}\int_0^{1-e^{-c}}\big(\frac{\log 2}{1-x}-\frac{\log(2-x)}{1-x}\big)dx=\frac1{c^2}\int_0^{1-e^{-c}}\frac{\log(1-\frac x2)}{x-1}dx,  \ q_n=1-\frac cn,\ c>0.
\end{aligned}
\end{equation}

Now we turn to $E_\infty^{s;\text{Geo}(1-q_n)}\mathcal{I}_n$.
The expectation $E_\infty^{s;\text{Geo}(1-q_n)}\mathcal{I}_n$ is given by \eqref{representexps}.
Of course,
\begin{equation}\label{firsttermsagain}
\frac {q_n}{1-q_n}(n-1)\sim \frac{n^2}c.
\end{equation}
From \eqref{asymptotically} and \eqref{evalint}, we have
\begin{equation}
\sum_{j=1}^{n-1}\frac{jq_n^j}{1-q_n^j}\sim\frac{n^2}{c^2}\int_0^c\frac{ye^{-y}}{1-e^{-y}}dy.
\end{equation}
By a change of variables, we have
\begin{equation}\label{changevar}
\int_0^c\frac{ye^{-y}}{1-e^{-y}}dy=-\int_0^{1-e^{-c}}\frac{\log(1-x)}xdx.
\end{equation}
From  \eqref{representexps} and \eqref{firsttermsagain}-\eqref{changevar}, we conclude
that
\begin{equation}\label{finalalpha1s}
\begin{aligned}
&\lim_{n\to\infty}\frac{E_\infty^{s;\text{Geo}(1-q_n)}\mathcal{I}_n}{n^2}=\frac1c+\frac1{c^2}\int_0^{1-e^{-c}}\frac{\log(1-x)}xdx=\\
&\frac1{c^2}\int_0^{1-e^{-c}}\Big(\frac1{1-x}+\frac{\log(1-x)}x\Big)dx.
\end{aligned}
\end{equation}

\section{Proofs of Proposition \ref{fixed} and Theorem \ref{1}}\label{wlln-sec}
\it\noindent  Proof of Proposition \ref{fixed}.\rm\
For the shifted case, we represent $\mathcal{I}_n$ as $\mathcal{I}=\sum_{j=2}^nI_{<j}$, where $I_{<j}$ is the number of inversions involving pairs $\{\{i,j\}:1\le i<j\}$.
In the shifted case,  by Proposition \ref{indepinver} and the remark following it, the random variables $\{1_{<j}\}_{j=2}^\infty$ are independent
and have
truncated binomial distributions with fixed parameter $1-q$; thus their variances are uniformly bounded. Denoting variance in the shifted case by
$\text{Var}_{s;1-q}$,
we have $\text{Var}_{s;1-q}(\mathcal{I}_n)=\sum_{j=2}^n\text{Var}_{s;1-q}(I_{<j})\le Cn$, for some constant $C$.
In section \ref{expectation-sec} we showed that  with fixed $q$, the expected value of $\mathcal{I}_n$ in the shifted
 case is on the order $n$.
Thus, by the second moment method,
\begin{equation}\label{s-fixed}
\text{w}-\lim_{n\to\infty}\frac{\mathcal{I}_n}{E^{s;\text{Geo}(1-q)}_\infty\mathcal{I}_n}=1 \ \text{under}\ P_\infty^{s;\text{Geo}(1-q)}.
\end{equation}
Proposition \ref{fixed} for the shifted case follows from \eqref{s-fixed} and \eqref{asympexp-s}.

Let $\text{Var}_{b;1-q}$ denote  variance in the biased case.
In section \ref{expectation-sec} we showed that  with fixed $q$, the expected value of $\mathcal{I}_n$ in the biased
 case is on the order $n$.
We will show that $\text{Var}_{b;1-q}(\mathcal{I}_n)$  is also on the order $n$.

It is clear from the biased construction that
$1_{\{\sigma^{-1}_j<\sigma^{-1}_i\}}$ and $1_{\{\sigma^{-1}_l<\sigma^{-1}_k\}}$ are independent if
$\{i,j\}\cap\{k,l\}=\emptyset$.
Writing $\mathcal{I}_n=\sum_{1\le i<j\le n}1_{\sigma^{-1}_j<\sigma^{-1}_i}$,
we have
\begin{equation*}
\begin{aligned}
&E_\infty^{b;\text{Geo}(1-q)}(\mathcal{I}_n)^2=\sum_{1\le i<j\le n}\sum_{1\le k<l\le n}E_\infty^{b;\text{Geo}(1-q)}1_{\{\sigma^{-1}_j<\sigma^{-1}_i\}}1_{\{\sigma^{-1}_l<\sigma^{-1}_k\}}=\\
&\sum_{1\le i<j\le n}\Big(\sum_{1\le k<l\le n:\{i,j\}\cap\{k,l\}=\emptyset}
E_\infty^{b;\text{Geo}(1-q)}1_{\{\sigma^{-1}_j<\sigma^{-1}_i\}}
E_\infty^{b;\text{Geo}(1-q)}1_{\{\sigma^{-1}_l<\sigma^{-1}_k\}}\Big)+\\
&\sum_{1\le i<j\le n}\Big(\sum_{1\le k<l\le n:\{i,j\}\cap\{k,l\}\neq\emptyset}
E_\infty^{b;\text{Geo}(1-q)}1_{\{\sigma^{-1}_j<\sigma^{-1}_i\}}1_{\{\sigma^{-1}_l<\sigma^{-1}_k\}}\Big)\le\\
&(E_\infty^{b;\text{Geo}(1-q)}\mathcal{I}_n)^2+
\sum_{1\le i<j\le n}\Big(\sum_{1\le k<l\le n:\{i,j\}\cap\{k,l\}\neq\emptyset}
E_\infty^{b;\text{Geo}(1-q)}1_{\{\sigma^{-1}_j<\sigma^{-1}_i\}}1_{\{\sigma^{-1}_l<\sigma^{-1}_k\}}\Big).
\end{aligned}
\end{equation*}
Thus,
\begin{equation}\label{varianceest}
\text{Var}_{b;1-q}(\mathcal{I}_n)\le \sum_{1\le i<j\le n}\Big(\sum_{1\le k<l\le n:\{i,j\}\cap\{k,l\}\neq\emptyset}
E_\infty^{b;\text{Geo}(1-q)}1_{\{\sigma^{-1}_j<\sigma^{-1}_i\}}1_{\{\sigma^{-1}_l<\sigma^{-1}_k\}}\Big).
\end{equation}
We break the sum on the right hand side of \eqref{varianceest} into five parts, depending on the values of $(k,l)$. The first part
is with $(k,l)$ satisfying $l=j$ and $k\neq i$; the second part is with $l=i$; the third
part is with $k=j$; the fourth part is with $k=i$ and $l\neq j$; and the fifth part
is with $(k,l)=(i,j)$.

The fifth part is equal to $E_\infty^{b;\text{Geo}(1-q)}\mathcal{I}_n$, so it is of order $n$. We will now show that
each of the first four parts is also of order $n$. Denote the $i$th part by $I_i(n)$.
For the first part, since $l=j$, we have $1\le k<j$ as well as $k\neq i$. Thus
$I_1(n)=\sum_{1\le i,k<j\le n; k\neq i}
E_\infty^{b;\text{Geo}(1-q)}1_{\{\sigma^{-1}_j<\sigma^{-1}_i\}}1_{\{\sigma^{-1}_j<\sigma^{-1}_k\}}$.
We have
$$
E_\infty^{b;\text{Geo}(1-q)}1_{\{\sigma^{-1}_j<\sigma^{-1}_i\}}1_{\{\sigma^{-1}_j<\sigma^{-1}_k\}}=\frac{p_j}{p_i+p_j+p_k}=
\frac{q^j}{q^i+q^j+q^k}.
$$
Therefore
\begin{equation}\label{I1}
I_1(n)\le\sum_{1\le i,k<j\le n}\frac{q^j}{q^i+q^j+q^k}.
\end{equation}
By Riemann sum approximation, we have
\begin{equation}\label{riemupper}
\begin{aligned}
&\sum_{1\le k<j}\frac{q^j}{q^i+q^j+q^k}\le\int_0^{j-1}\frac{q^j}{q^i+q^j+e^{x\log q}}dx\le\\
&\int_0^j\frac{q^je^{-x\log q}}{1+(q^i+q^j)e^{-x\log q}}dx\le\frac{q^j}{(-\log q)(q^j+q^i)}\log(2+q^{i-j}).
\end{aligned}
\end{equation}
From \eqref{I1} and \eqref{riemupper} we have
\begin{equation}
\begin{aligned}
&I_1(n)\le\frac1{-\log q} \sum_{1\le i<j\le n}\frac{q^j}{(q^j+q^i)}\log(2+q^{i-j})=\\
&\frac1{-\log q}\sum_{r=1}^{n-1}(n-r)\frac{q^r}{1+q^r}\log(2+q^{-r})\le\frac n{-\log q}\sum_{r=1}^{n-1}
\frac{q^r}{1+q^r}\big(C+(-\log q)r\big)\le C_1n,
\end{aligned}
\end{equation}
for constants $C,C_1>0$.

The other three parts follow similarly.
Indeed
$$
\begin{aligned}
&I_2(n)=\sum_{1\le k< i<j\le n}E_\infty^{b;\text{Geo}(1-q)}1_{\{\sigma^{-1}_j<\sigma^{-1}_i\}}
1_{\{\sigma^{-1}_i<\sigma^{-1}_k\}}=\\
&\sum_{1\le k< i<j\le n}\frac{p_j}{p_i+p_j+p_k}\frac{p_i}{p_i+p_k}\le
\sum_{1\le k< i<j\le n}\frac{q^j}{q^i+q^j+q^k},
\end{aligned}
$$
and the right hand side above is less than the right hand side of \eqref{I1}.
Also,
$$
\begin{aligned}
&I_3(n)=\sum_{1\le i<j<l\le n}E_\infty^{b;\text{Geo}(1-q)}1_{\{\sigma^{-1}_j<\sigma^{-1}_i\}}
1_{\{\sigma^{-1}_l<\sigma^{-1}_j\}}=\\
&\sum_{1\le i<j<l\le n}\frac{p_l}{p_i+p_j+p_l}\frac{p_j}{p_i+p_j}\le
\sum_{1\le i<j<l\le n}\frac{q^l}{q^i+q^j+q^l},
\end{aligned}
$$
and the right hand side above is less than the  right hand side of \eqref{I1}.
Finally,
$$
\begin{aligned}
&I_4(n)=\sum_{1\le i<j\le n, \thinspace l\in\{i+1,\cdots,n\}-\{j\}}
E_\infty^{b;\text{Geo}(1-q)}1_{\{\sigma^{-1}_j<\sigma^{-1}_i\}}
1_{\{\sigma^{-1}_l<\sigma^{-1}_i\}}=\\
&\sum_{1\le i<j\le n, \thinspace l\in\{i+1,\cdots,n\}-\{j\}}\big(\frac{p_j}{p_i+p_j+p_l}\frac{p_l}{p_i+p_l}
+\frac{p_l}{p_i+p_j+p_l}\frac{p_j}{p_i+p_j}\big)\le\\
&2\sum_{1\le i<j,l\le n}\frac{q^j}{q^i+q^j}\frac{q^l}{q^i+q^l}=2\sum_{1\le i<j\le n}
\frac{q^j}{q^i+q^j}\sum_{i<l\le n}\frac{q^l}{q^i+q^l}=\\
&2\sum_{1\le i<j\le n}
\frac{q^j}{q^i+q^j}\sum_{r=1}^{n-i}\frac{q^r}{1+q^r}\le
C\sum_{1\le i<j\le n}
\frac{q^j}{q^i+q^j}=
CE_\infty^{b;\text{Geo}(1-q)}\mathcal{I}_n,
\end{aligned}
$$
for some $C>0$.

Since $\text{Var}_{b;1-q}(\mathcal{I}_n)$ is on the order $n$, by the second moment method,
\begin{equation}\label{b-fixed}
\text{w}-\lim_{n\to\infty}\frac{\mathcal{I}_n}{E^{b;\text{Geo}(1-q)}_\infty}=1 \ \text{under}\ P_\infty^{b;\text{Geo}(1-q)}.
\end{equation}
Proposition \ref{fixed} for the biased case then follows from \eqref{b-fixed} along with \eqref{asympexp-b} and \eqref{q-b}.
\hfill $\square$
\medskip

\noindent\it Proof of Theorem \ref{1}.\rm\ Consider $q_n$ as in part (a) or part (b).
For the shifted case, we use the same method of  proof  used for the shifted case  in Proposition \ref{fixed}.
Let $\text{Var}_{s;1-q_n}$ denote  variance in the shifted case.
We represent $\mathcal{I}_n$ as $\mathcal{I}=\sum_{j=2}^nI_{<j}$, where $I_{<j}$ is the number of inversions involving pairs $\{\{i,j\}:1\le i<j\}$.
By Proposition \ref{indepinver} and the remark following it, the random variables $\{1_{<j}\}_{j=2}^\infty$ are independent
and have
truncated binomial distributions with  parameter $1-q_n$. Thus, under the assumption of part (a),
$\text{Var}_{s;1-q_n}(1_{<j})\le Cn^{2\alpha}$, for some $C>0$ and all $j$, while under the assumption of part (b) the same inequality holds with $\alpha=1$.
Consequently,  $\text{Var}_{s;1-q_n}(\mathcal{I}_n)\le Cn^{1+2\alpha}$ under the assumption of part (a), while under the assumption of part (b)
  the same inequality holds with $\alpha=1$.
  In section \ref{expectation-sec} we showed that $E^{s;\text{Geo}(1-q_n)}\mathcal{I}_n$ is on the order $n^{1+\alpha}$ under the assumption of part (a), and on the order
  $n^2$ under the assumption of part (b).
  Therefore, both in parts (a) and (b) we have $\text{Var}_{s;1-q_n}(\mathcal{I}_n)=o\big((E_\infty^{s;\text{Geo}(1-q_n)}\mathcal{I}_n)^2\big)$.
Thus, by the second moment method,
\begin{equation}\label{s-notfixed}
\text{w}-\lim_{n\to\infty}\frac{\mathcal{I}_n}{E^{s;\text{Geo}(1-q_n)}_\infty}=1 \ \text{under}\ P_\infty^{s;\text{Geo}(1-q_n)}.
\end{equation}
The weak law stated in part (a)  for the shifted case follows from \eqref{s-notfixed} along with \eqref{finalalphas}, while the weak law stated in part (b) for the shifted case follows from
\eqref{s-notfixed} and \eqref{finalalpha1s}.

Now consider the biased case. Let $\text{Var}_{b;1-q_n}$ denote  variance in the biased case.
In  the biased case, it is clear from the construction that
$1_{\{\sigma^{-1}_j<\sigma^{-1}_i\}}$ and $1_{\{\sigma^{-1}_l<\sigma^{-1}_k\}}$ are independent if
$\{i,j\}\cap\{k,l\}=\emptyset$.
 Writing $\mathcal{I}_n=\sum_{1\le i<j\le n}1_{\sigma^{-1}_j<\sigma^{-1}_i}$,
we have
\begin{equation}
\begin{aligned}
&E_\infty^{b;q_n}(\mathcal{I}_n)^2=\sum_{1\le i<j\le n}\sum_{1\le k<l\le n}E_\infty^{b;q_n}1_{\{\sigma^{-1}_j<\sigma^{-1}_i\}}1_{\{\sigma^{-1}_l<\sigma^{-1}_k\}}=\\
&\sum_{1\le i<j\le n}\Big(\sum_{1\le k<l\le n:\{i,j\}\cap\{k,l\}=\emptyset}E_\infty^{b;q_n}1_{\{\sigma^{-1}_j<\sigma^{-1}_i\}}E_\infty^{b;q_n}1_{\{\sigma^{-1}_l<\sigma^{-1}_k\}}\Big)+\\
&\sum_{1\le i<j\le n}\Big(\sum_{1\le k<l\le n:\{i,j\}\cap\{k,l\}\neq\emptyset}E_\infty^{b;q_n}1_{\{\sigma^{-1}_j<\sigma^{-1}_i\}}1_{\{\sigma^{-1}_l<\sigma^{-1}_k\}}\Big)\le\\
&(E_\infty^{b;q_n}\mathcal{I}_n)^2+4n\sum_{1\le i<j\le n}E_\infty^{b;q_n}1_{\{\sigma^{-1}_j<\sigma^{-1}_i\}}=(E_\infty^{b;q_n}\mathcal{I}_n)^2+4nE_\infty^{b;q_n}\mathcal{I}_n.
\end{aligned}
\end{equation}
Thus $\text{Var}_{b;1-q_n}(\mathcal{I}_n)=O(nE_\infty^{b;q_n}\mathcal{I}_n)$.
In the cases of $q_n$ as in  parts (a) and (b) of the theorem, $E_\infty^{b;q_n}\mathcal{I}_n$
is on a larger order than $n$. Consequently, it follows that $\text{Var}_{b;1-q_n}(\mathcal{I}_n)=o\big((E_\infty^{b;q_n}\mathcal{I}_n)^2\big)$.
Thus, by the second moment method,
\eqref{s-notfixed} holds with $s$ replaced by $b$.
Using this with \eqref{finalnalpha}  proves the weak law stated in part (a) for the biased case, while using this with
\eqref{finalnalpha1s}  proves the weak law stated in  part (b) for the biased case.

This completes the proof of part (a),  and it completes the proof of part (b) except for the statement concerning the behavior of $I_b(c)$ and $I_s(c)$.
We leave it to the reader to check the claim regarding the behavior of these two functions as $c\to0$ and as $c\to\infty$.
It remains to show that
$I_b(c)<I_s(c)$. Of course, $I_b(c)\le I_s(c)$ follows by the stochastic dominance  in Proposition \ref{dom}.
It suffices to show that
$$
\frac1{1-x}+\frac{\log(1-x)}x+\frac{\log(1-\frac x2)}{1-x}>0,\ 0<x<1.
$$
Multiplying by $x(1-x)$, it suffices to show that
$$
F(x):=x+(1-x)\log(1-x)+x\log(1-\frac x2)>0,\ 0<x<1.
$$
We have $F(0)=0$. Differentiating gives
$$
F'(x)=-\log(1-x)+\log(1-\frac x2)-\frac x{2-x}.
$$
We have $F'(0)=0$.
Differentiating again gives
$$
F''(x)=\frac1{1-x}-\frac2{2-x}-\frac x{(2-x)^2}.
$$
We have $F''(0)=0$.
Differentiating a third time gives
$$
F'''(x)=\frac1{(1-x)^2}-\frac3{(2-x)^2}-\frac{2x}{(2-x)^3}=\frac{2+x-2x^2}{(1-x)^2(2-x)^3}>0,\ 0<x<1.
$$
This completes the proof of part (b).

We now turn to part (c).
For $q_1<q_2$ and $i<j$, it is immediate from the construction in the biased case and easy to check in the shifted case
(similar to the proof of Proposition \ref{dom}) that  $1_{\sigma^{-1}_j<\sigma^{-1}_i}$ under
$P_\infty^{*;\text{Geo}(1-q_2)}$ strictly stochastically dominates
$1_{\sigma^{-1}_j<\sigma^{-1}_i}$ under
$P_\infty^{*;\text{Geo}(1-q_1)}$, for $*=b$ or $*=s$.
 Thus, by the linearity of the expectation,
   $E_\infty^{*;\text{Geo}(1-q_n)}\mathcal{I}_n$  is smaller for
 $q_n=1-\frac cn$ with $c>0$ than it is for $q_n$ as in  part (c) and $n$ sufficiently large, where $*=b$ or $*=s$.
  By part (b),
 $\lim_{n\to\infty}E_\infty^{b;\text{Geo}(1-q_n)}\frac{E\mathcal{I}_n}{n^2}=I_b(c)$
and $\lim_{n\to\infty}E_\infty^{s;\text{Geo}(1-q_n)}\frac{E\mathcal{I}_n}{n^2}=I_s(c)$, and
  $\lim_{c\to0}I_b(c)=\lim_{c\to0}I_s(c)=\frac14$.
Thus, $\limsup_{n\to\infty}\frac{E_\infty^{*;\text{Geo}(1-q_n)}\mathcal{I}_n}{n^2}\le\frac14$, for
$q_n$ as in part (c) and $*=b$ or $*=s$.
On the other hand, a uniformly random permutation of $S_n$ can also be constructed via the biased
or shifted constructions, by letting $p_j=\frac1n$, $j=1,\cdots, n$, and then the same consideration
as in the first line of this paragraph shows that
the expected value $E\mathcal{I}_n$ in the uniform case is larger than
$E_\infty^{*;\text{Geo}(1-q_n)}\mathcal{I}_n$, for $q_n$ in part (c) and $*=b$ or $*=s$.
It is well-known that in the uniform distribution case, $\lim_{n\to\infty}\frac{E\mathcal{I}_n}{n^2}=\frac14$.
 Part (c) follows from the above considerations.
 \hfill $\square$

\section{Proof of Theorem \ref{extra}}\label{extrastuff}
Let $\{w_k\}_{k=1}^\infty$ be a realization of the IID Beta$(1,\alpha)$-distributed random variables $\{W_k\}_{k=1}^\infty$, and let $\{p_k\}_{k=1}^\infty$ denote the corresponding realization of $\{\mathcal{P}_k\}_{k=1}^\infty$. So
\begin{equation}\label{thepk}
p_k=w_k\prod_{i=1}^{k-1}(1-w_i),\ k=1,2,\cdots.
\end{equation}
By Proposition \ref{indepinver}, under
$P_\infty^{s;\{p_k\}}$,
the random variables $\{1_{<j}\}_{j=2}^\infty$ are independent  and distributed according
to \eqref{<jdist}. In particular then, under $P_\infty^{s;\{p_k\}}$ these random variables converge in distribution as $j\to\infty$ to a
random variable $X$ with distribution
$P(X=k)=p_{k+1},\ k=0,1,\cdots$.
From \eqref{thepk}, we write
$p_k=w_ke^{\sum_{i=1}^{k-1}\log(1-w_i)}$ and note that by the law of large numbers,
$\frac1k\sum_{i=1}^k\log(1-w_i)$ converges $P_\theta$-almost surely as $k\to\infty$ to $E_\theta\log (1-W_1)<0$.
Consequently, $P_\theta$-almost surely, the $\{p_k\}_{k=1}^\infty$
 decay exponentially.
Therefore, $EX^2<\infty$ $P_\theta$-almost surely.
Since the distributions of the $\{1_{<j}\}_{j=2}^\infty$ are truncated versions of the distribution of $X$,
the random variable $X$ stochastically dominates all of the $\{1_{<j}\}_{j=2}^\infty$. Thus,
the second moments of the $\{1_{<j}\}_{j=2}^\infty$ are $P_\theta$-almost surely uniformly bounded.
We have
$\lim_{j\to\infty}E_\infty^{s;\{p_k\}}1_{<j}=E_\theta X=\sum_{k=1}^\infty kp_k$,
 $P_\theta$-almost surely.
From these facts, we conclude that $P_\theta$-almost surely, the weak law of large numbers holds
for $\{1_{<j}\}_{j=2}^\infty$ in the form
$w-\lim_{n\to\infty}\frac1n\sum_{j=2}^n1_{<j}=EX=\sum_{k=1}^\infty kp_{k+1}$.
Using this with
 \eqref{thepk} and the fact that
$\mathcal{I}_n=\sum_{j=2}^n1_{<j}$, we obtain \eqref{weakgemshift}.

We now prove \eqref{expextra}.
From the previous paragraph and \eqref{<jdist}, we have
\begin{equation}\label{propgemshift1}
E_\infty^{s;\text{GEM}(\theta)}1_{<j}=E_\theta\sum_{k=1}^{j-1}k\frac{W_{k+1}\prod_{i=1}^k(1-W_i)}{\sum_{k=1}^j\mathcal{P}_k}.
\end{equation}
Also,
\begin{equation}\label{propgemshift2}
\lim_{j\to\infty}\sum_{k=1}^{j-1}k\frac{W_{k+1}\prod_{i=1}^k(1-W_i)}{\sum_{k=1}^j\mathcal{P}_k}=
\sum_{k=1}^\infty kW_{k+1}\prod_{i=1}^k(1-W_i), \ P_\theta-\text{almost surely}.
\end{equation}
Recalling that $\mathcal{P}_1=W_1$ and $\mathcal{P}_2=(1-W_1)W_2$, we have
\begin{equation}\label{propgemshift3}
\sum_{k=1}^{j-1}k\frac{W_{k+1}\prod_{i=1}^k(1-W_i)}{\sum_{k=1}^j\mathcal{P}_k}\le
\sum_{k=1}^\infty k\frac{W_{k+1}\prod_{i=1}^k(1-W_i)}{W_1+(1-W_1)W_2},\ \text{for all}\ j\ge2.
\end{equation}
We will show that
\begin{equation}\label{dct}
E_\theta\sum_{k=1}^\infty k\frac{W_{k+1}\prod_{i=1}^k(1-W_i)}{W_1+(1-W_1)W_2}<\infty.
\end{equation}
It then follows from \eqref{propgemshift1}-\eqref{dct} and the dominated convergence theorem that
\begin{equation}\label{dctagain}
\lim_{j\to\infty}E_\infty^{s;\text{GEM}(\theta)}1_{<j}=\sum_{k=1}^\infty kE_\theta W_{k+1}\prod_{i=1}^k(1-W_i).
\end{equation}
A straightforward calculation will reveal that
\begin{equation}\label{propgemshift4}
\sum_{k=1}^\infty kE_\theta W_{k+1}\prod_{i=1}^k(1-W_i)=\theta.
\end{equation}
Since $E_\infty^{s;\text{GEM}(\theta)}\mathcal{I}_n=\sum_{j=2}^n E_\infty^{s;\text{GEM}(\theta)}1_{<j}$,
it then follows from \eqref{dctagain} and \eqref{propgemshift4} that
$\lim_{n\to\infty}\frac1nE_\infty^{s;\text{GEM}(\theta)}\mathcal{I}_n=\theta$, completing the proof of \eqref{expextra}.
Thus, it remains to prove \eqref{dct} and \eqref{propgemshift4}.

We have
$$
E_\theta W_{k+1}\prod_{i=1}^k(1-W_i)=E_\theta W_1\big(E_\theta(1-W_1)\big)^k=\frac1{1+\theta}\big(\frac\theta{1+\theta})^k.
$$
Thus,
$$
\begin{aligned}
&\sum_{k=1}^\infty kE_\theta W_{k+1}\prod_{i=1}^k(1-W_i)=\frac1{1+\theta}\sum_{k=1}^\infty k\big(\frac\theta{1+\theta})^k=\\
&\frac\theta{(1+\theta)^2}\frac{d}{d\lambda}(\frac1{1-\lambda})|_{\lambda=\frac\theta{1+\theta}}=\theta,
\end{aligned}
$$
proving \eqref{propgemshift4}.

We now turn to \eqref{dct}.
For the $k$th summand in \eqref{dct}, we have
\begin{equation}\label{summands3}
E_\theta\frac{W_{k+1}\prod_{i=1}^k(1-W_i)}{W_1+(1-W_1)W_2}=
E_\theta\frac{(1-W_1)(1-W_2)}{W_1+(1-W_1)W_2}E_\theta W_{k+1}\prod_{i=3}^k(1-W_i),\ \text{for}\ k\ge 3,
\end{equation}
while for $k=2$ we have
\begin{equation}\label{summand2}
E_\theta\frac{W_2(1-W_1)}{W_1+(1-W_1)W_2}\le 1.
\end{equation}
We have
\begin{equation}\label{summands3part2}
E_\theta W_{k+1}\prod_{i=3}^k(1-W_i)=E_\theta W_1\big(E_\theta (1-W_1)\big)^{k-2}=
\frac1{1+\theta}\big(\frac\theta{1+\theta}\big)^{k-2}.
\end{equation}
And finally,
\begin{equation}\label{summands3part1}
\begin{aligned}
&E_\theta\frac{(1-W_1)(1-W_2)}{W_1+(1-W_1)W_2}=\theta^2
\int_0^1dw_1\int_0^1dw_2\frac{(1-w_1)^\theta(1-w_2)^\theta}{w_1+(1-w_1)w_2}\le\\
&\theta^2\int_0^1dw_1\int_0^1dw_2\frac{(1-w_1)^\theta}{w_1+(1-w_1)w_2}=\\
&\theta^2\int_0^1dw_1(1-w_1)^{\theta-1}\log\big(w_1+(1-w_1)w_2\big)|_{w_2=0}^1=\\
&-\theta^2\int_0^1(1-w_1)^{\theta-1}\log w_1\thinspace dw_1<\infty.
\end{aligned}
\end{equation}
Now \eqref{dct} follows from \eqref{summands3}-\eqref{summands3part1}.
\hfill $\square$

\section{Proof of Theorem \ref{biasedgemthm}}\label{GEMofaproof}
Recall that  $P_\theta$ and $E_\theta$  denote respectively probability and   expectation with respect to the IID Beta$(1,\theta)$-distributed sequence $\{W_k\}_{k=1}^\infty$ that is associated with the GEM$(\theta)$ distribution.
Analogous to the first paragraph of section \ref{expectation-sec}, to calculate the expected number of inversions, we write $\mathcal{I}_n=\sum_{1\le i<j\le n}1_{\{\sigma^{-1}_j<\sigma^{-1}_i\}}$. It is immediate from the construction that
\begin{equation}\label{basicformula}
\begin{aligned}
&E_\infty^{b;\text{\rm GEM}(\theta)}1_{\{\sigma^{-1}_j<\sigma^{-1}_i\}}=E_\theta\frac{(1-W_1)\cdots(1-W_{j-1})W_j}{(1-W_1)\cdots(1-W_{i-1})W_i+(1-W_1)\cdots(1-W_{j-1})W_j}=\\
&1-E_\theta\frac1{1+\frac{1-W_i}{W_i}(1-W_{i+1})\cdots(1-W_{j-1})W_j}=\\
&1-E_\theta\frac1{1+\frac{1-W_1}{W_1}(1-W_2)\cdots(1-W_k)W_{k+1}},\ \ k=j-i.
\end{aligned}
\end{equation}
Thus,
\begin{equation}\label{expfirstform}
E_\infty^{b;\text{\rm GEM}(\theta)}\mathcal{I}_n=\sum_{k=1}^{n-1}(n-k)\big(1-E_\theta\frac1{1+\frac{1-W_1}{W_1}(1-W_2)\cdots(1-W_k)W_{k+1}}\big).
\end{equation}
We will show that
\begin{equation}\label{needtoshow}
\sum_{k=1}^\infty\big(1-E_\theta\frac1{1+\frac{1-W_1}{W_1}(1-W_2)\cdots(1-W_k)W_{k+1}}\big)=\theta.
\end{equation}
From \eqref{expfirstform} and \eqref{needtoshow} it follows that
$$
\lim_{n\to\infty}\frac{E_\infty^{b;\text{\rm GEM}(\theta)}\mathcal{I}_n}n=\theta.
$$
Indeed, note
that the summands in \eqref{needtoshow} are positive, which follows from \eqref{basicformula}, and note
 from \eqref{needtoshow} that for any $\epsilon>0$, there exists a $K_\epsilon$ such that
\newline
$\sum_{k=k_0}^\infty\big(1-E_\theta\frac1{1+\frac{1-W_1}{W_1}(1-W_2)\cdots(1-W_k)W_{k+1}}\big)<\epsilon,\ \text{for}\ k_0> K_\epsilon$.
Thus, for\newline $n>K_\epsilon$,
$$
\begin{aligned}
&\sum_{k=1}^{n-1}k\big(1-E_\theta\frac1{1+\frac{1-W_1}{W_1}(1-W_2)\cdots(1-W_k)W_{k+1}}\big)\le K_\epsilon \theta+\epsilon n.
\end{aligned}
$$

To complete the proof of the theorem,
we now turn to  the proof of \eqref{needtoshow}. We calculate the density $f_{\frac{1-W_1}{W_1}}(z)$ of the random variable $\frac{1-W_1}{W_1}$.
We have
$$
P_\theta(\frac{1-W_1}{W_1}\le z)=P_\theta(W_1\ge\frac1{1+z})=\int_{(1+z)^{-1}}^1\theta(1-w)^{\theta-1}dw,
$$
from which it follows that
$$
f_{\frac{1-W_1}{W_1}}(z)=\frac{\theta z^{\theta-1}}{(1+z)^{1+\theta}},\ 0<z<\infty.
$$
Letting
$$
\alpha_k=(1-W_2)\cdots(1-W_k)W_{k+1},\ k\ge1,
$$
we have
\begin{equation}\label{firstint}
E_\theta\frac1{1+\frac{1-W_1}{W_1}(1-W_2)\cdots(1-W_k)W_{k+1}}=\theta E_\theta\int_0^\infty\frac1{1+\alpha_kz}\thinspace \frac{z^{\theta-1}}{(1+z)^{\theta+1}}dz.
\end{equation}
Making the substitution $u=\frac z{1+z}$, we obtain
\begin{equation}\label{subst}
\theta\int_0^\infty\frac1{1+\alpha_kz}\thinspace \frac{z^{\theta-1}}{(1+z)^{\theta+1}}dz=\theta\int_0^1\frac{u^{\theta-1}(1-u)}{1-u+\alpha_k u}du=
1-\theta\alpha_k\int_0^1\frac{u^\theta}{1-u+\alpha u}du.
\end{equation}
From \eqref{firstint} and \eqref{subst} we have
\begin{equation}\label{stage1}
1-E_\theta\frac1{1+\frac{1-W_1}{W_1}(1-W_2)\cdots(1-W_k)W_{k+1}}=\theta E_\theta\alpha_k\int_0^1\frac{u^\theta}{1-u+\alpha_k u}du.
\end{equation}

We now write
\begin{equation}\label{stage2}
\begin{aligned}
&E_\theta\alpha_k\int_0^1\frac{u^\theta}{1-u+\alpha_k u}du=E_\theta\alpha_k\int_0^1u^\theta\Big(\sum_{m=0}^\infty u^m(1-\alpha_k)^m\Big)du=\\
&E_\theta\sum_{m=0}^\infty\frac{\alpha_k}{m+\theta+1}(1-\alpha_k)^m=\sum_{m=0}^\infty\frac1{m+\theta+1}\Big(\sum_{i=0}^m(-1)^i\binom miE_\theta\alpha_k^{i+1}\Big).
\end{aligned}
\end{equation}
We have
\begin{equation}\label{3exps-1}
E_\theta\alpha_k^{i+1}=\big(E_\theta(1-W_1)^{i+1}\big)^{k-1}E_\theta W_1^{i+1}.
\end{equation}
Also,
\begin{equation}\label{3exps-2}
E_\theta(1-W_1)^{i+1}=\int_0^1(1-w)^{i+1}\theta(1-w)^{\theta-1}dw=\frac\theta{\theta+i+1},
\end{equation}
and from the well-known normalization for the Beta-distributions,
\begin{equation}\label{3exps-3}
\begin{aligned}
&E_\theta W_1^{i+1}=\int_0^1w^{i+1}\theta(1-w)^{\theta-1}dw=\frac{\theta\Gamma(\theta)\Gamma(i+2)}{\Gamma(\theta+i+2)}=\\
&\frac{\Gamma(\theta+1)(i+1)!}{\Gamma(\theta+i+2)}=\frac{(i+1)!}{\prod_{l=1}^{i+1}(\theta+l)}=\frac1{\binom {\theta+i+1}{i+1}}.
\end{aligned}
\end{equation}
Substituting \eqref{3exps-1}-\eqref{3exps-3}   in \eqref{stage2}, and using this with \eqref{stage1}, we obtain
\begin{equation}\label{stage3}
\begin{aligned}
&1-E_\theta\frac1{1+\frac{1-W_1}{W_1}(1-W_2)\cdots(1-W_k)W_{k+1}}=\\
&\theta\sum_{m=0}^\infty\frac1{m+\theta+1}\Big(\sum_{i=0}^m(-1)^i\frac{\binom mi}{\binom{\theta+i+1}{i+1}}\big(\frac{\theta}{\theta+i+1}\big)^{k-1}\Big).
\end{aligned}
\end{equation}

Recall from the above calculations that
$\sum_{i=0}^m(-1)^i\frac{\binom mi}{\binom{\theta+i+1}{i+1}}\big(\frac{\theta}{\theta+i+1}\big)^{k-1}=E_\theta\alpha_k(1-\alpha_k)^m>0$.
Thus, summing \eqref{stage3}
over $k$ and invoking the monotone convergence theorem, we obtain
\begin{equation}\label{stage4}
\begin{aligned}
&\sum_{k=1}^\infty\big(1-E_\theta\frac1{1+\frac{1-W_1}{W_1}(1-W_2)\cdots(1-W_k)W_{k+1}}\big)=\\
&\theta\sum_{m=0}^\infty\frac1{m+\theta+1}\Big(\sum_{i=0}^m(-1)^i\frac{\binom mi}{\binom{\theta+i+1}{i+1}}\frac{\theta+i+1}{i+1}\Big)=\\
&\theta\sum_{m=0}^\infty\frac1{m+\theta+1}\Big(\sum_{i=0}^m(-1)^i\frac{\binom mi}{\binom{\theta+i}i}\Big).
\end{aligned}
\end{equation}
In light of \eqref{stage4}, to complete the proof of
\eqref{needtoshow} we need to show that
\begin{equation}\label{equalto1}
\sum_{m=0}^\infty\frac1{m+\theta+1}\Big(\sum_{i=0}^m(-1)^i\frac{\binom mi}{\binom{\theta+i}i}\Big)=1,\ \theta>0.
\end{equation}

We first prove \eqref{equalto1} for $\theta\in\mathbb{N}$.
When $\theta\in\mathbb{N}$, we can write
$\binom{\theta+i}i=\binom{\theta+i}\theta=\frac{(\theta+i)!}{\theta!\thinspace i!}$. Thus,
\begin{equation}\label{stage5}
\begin{aligned}
&\sum_{i=0}^m(-1)^i\frac{\binom mi}{\binom{\theta+i}i}=\theta!
\sum_{i=0}^m(-1)^i\frac{m!}{(m-i)!}\frac1{(\theta+i)!}=\\
&\frac{\theta!}{(m+1)\cdots(m+\theta)}\sum_{i=0}^m(-1)^i\binom{m+\theta}{\theta+i}=\\
&(-1)^{\theta-1}\frac{\theta!}{(m+1)\cdots(m+\theta)}\sum_{j=0}^{\theta-1}(-1)^j\binom{m+\theta}j,\ \ \theta\in\mathbb{N},
\end{aligned}
\end{equation}
where the last equality follows from the fact that $\sum_{j=0}^{m+\theta}(-1)^j\binom{m+\theta}j=0$.

We now show that
\begin{equation}\label{equals1}
(-1)^{\theta-1}\frac{(\theta-1)!}{(m+1)\cdots(m+\theta-1)}\sum_{j=0}^{\theta-1}(-1)^j\binom{m+\theta}j=1.
\end{equation}
 Let
$$
f(m)=(-1)^{\theta-1}(\theta-1)!\sum_{j=0}^{\theta-1}(-1)^j\binom{m+\theta}j;\ \  \
g(m)=(m+1)\cdots(m+\theta-1).
$$
Both $f$ and $g$ are polynomials of degree $\theta-1$. They both have leading order coefficient equal to 1. The roots of $g$    are $\{-\theta+l\}_{l=1}^{\theta-1}$.
We now show that $f$ has the same roots, from which  \eqref{equals1} follows.
Of course it suffices to show that $h(m):=\sum_{j=0}^{\theta-1}(-1)^j\binom{m+\theta}j$ has the same roots.
We have
$$
h(-\theta+l)=\sum_{j=0}^{\theta-1}(-1)^j\binom lj=\sum_{j=0}^l(-1)^j\binom lj=0,\ l=1,\cdots \theta-1,
$$
where the   second equality follows from the fact that $\binom lj=0$, for \newline $j=l+1,\cdots, \theta-1$.

From \eqref{stage5} and \eqref{equals1} we have
\begin{equation}\label{final}
\sum_{i=0}^m(-1)^i\frac{\binom mi}{\binom{\theta+i}i}=\frac\theta{m+\theta},\ \ \theta\in\mathbb{N},\ m=0,1,\cdots.
\end{equation}
From \eqref{final} we conclude that
\begin{equation}\label{finalinteger}V
\begin{aligned}
&\sum_{m=0}^\infty\frac1{m+\theta+1}\Big(\sum_{i=0}^m(-1)^i\frac{\binom mi}{\binom{\theta+i}i}\Big)=\\
&\theta\sum_{m=0}^\infty\frac1{(m+\theta)(m+\theta+1)}=\theta
\sum_{m=0}^\infty\big(\frac1{m+\theta}-\frac1{m+\theta+1}\big)=1,\ \  \theta\in\mathbb{N}.
\end{aligned}
\end{equation}

We now show that \eqref{equalto1} in fact holds for all $\theta>0$. From \eqref{final} and \eqref{finalinteger},
it suffices to show that \eqref{final} holds for all $\theta>0$. Fix $m\in\{0,1,\cdots\}$.
Define
$$
A(\theta)=\sum_{i=0}^m(-1)^i\frac{\binom mi}{\binom{\theta+i}i}\thinspace; \ \ B(\theta)=\frac{\theta}{m+\theta}.
$$
Then $A$ is analytic
for $\theta\in\mathbb{C}-\{-l\}_{l=1}^m$, and $B$ is analytic for  $\theta\in\mathbb{C}-\{-m\}$.
Define $\mathcal{A}(\theta)=A(\frac1\theta)$ and $\mathcal{B}(\theta)=B(\frac1\theta)$.
Since $\lim_{\theta\to0}\mathcal{A}(\theta)=\lim_{\theta\to0}\mathcal{B}(\theta)=1$,
it follows that  $\theta=0$  is a removable singularity for $\mathcal{A}$ and $\mathcal{B}$. Hence, defining
$\mathcal{A}(0)=\mathcal{B}(0)=1$ makes $\mathcal{A}$ and $\mathcal{B}$ analytic functions
in a neighborhood of the origin.
Since $\mathcal{A}$ and $\mathcal{B}$ coincide on $\{0\}\cup\{\frac1n\}_{n=1}^\infty$, it follows from
the uniqueness theorem for analytic functions that $\mathcal{A}\equiv\mathcal{B}$ on
$\mathbb{C}-\{-l\}_{l=1}^m$, and thus in particular, $A(\theta)=B(\theta)$, for $\theta>0$.
\hfill $\square$

\section{Proof of Theorem \ref{indeptheta}}\label{proofindeptheta}
Let  the generic $P$ and $E$  denote respectively probability and   expectation with respect to the IID sequence $\{U_k\}_{k=1}^\infty$ of uniformly distributed random variables on $[0,1]$.
Analogous to the first paragraph of section \ref{expectation-sec}, to calculate the expected number of inversions, we write $\mathcal{I}_n=\sum_{1\le i<j\le n}1_{\{\sigma^{-1}_j<\sigma^{-1}_i\}}$. It is immediate from the construction that
\begin{equation}\label{basicformulaindep}
\begin{aligned}
&E_\infty^{b;\text{\rm IID-prod}(\theta)}1_{\{\sigma^{-1}_j<\sigma^{-1}_i\}}=
E\frac{\prod_{l=1}^jU_l^\frac1\theta}{\prod_{l=1}^iU_l^\frac1\theta+\prod_{l=1}^jU_l^\frac1\theta}=\\
&1-E\frac1{1+\prod_{l=i+1}^jU_l^\frac1\theta}=
1-E\frac1{1+\prod_{l=1}^kU_l^\frac1\theta},\ \ k=j-i.
\end{aligned}
\end{equation}
Thus,
\begin{equation}\label{firstformulaindep}
E_\infty^{b;\text{IID-prod}(\theta)}\mathcal{I}_n=\sum_{k=1}^n(n-k)\big(1-E\frac1{1+\prod_{l=1}^kU_l^\frac1\theta}\big).
\end{equation}
We will show that
\begin{equation}\label{toshowindep}
\sum_{k=1}^\infty\big(1-E\frac1{1+\prod_{l=1}^kU_l^\frac1\theta}\big)=\theta\log 2.
\end{equation}
Just as the displayed equation after \eqref{needtoshow} follows from  \eqref{expfirstform} and \eqref{needtoshow},
it follows from \eqref{firstformulaindep} and \eqref{toshowindep} that
\begin{equation}
\lim_{n\to\infty}\frac{E_\infty^{b;\text{\rm IID-prod}(r\theta)}\mathcal{I}_n}n=\theta\log2.
\end{equation}

To complete the proof of the theorem, we turn to the proof of \eqref{toshowindep}.
We have
\begin{equation}\label{calcindep}
\begin{aligned}
&E\frac1{1+\prod_{l=1}^kU_l^\frac1\theta}=E\sum_{m=0}^\infty(-1)^m\big(\prod_{l=1}^kU_l^\frac1\theta\big)^m=
\sum_{m=0}^\infty(-1)^m(EU_1^\frac m\theta)^k=\\
&\sum_{m=0}^\infty(-1)^m(\frac\theta{m+\theta})^k=1-\sum_{m=1}^\infty(-1)^{m-1}(\frac\theta{m+\theta})^k.
\end{aligned}
\end{equation}
From \eqref{calcindep} we have
\begin{equation}\label{ontheway}
\begin{aligned}
&\sum_{k=1}^\infty\big(1-E\frac1{1+\prod_{l=1}^kU_l^\frac1\theta}\big)=
\sum_{k=1}^\infty\sum_{m=1}^\infty(-1)^{m-1}(\frac\theta{m+\theta})^k=\\
&\lim_{K\to\infty}\lim_{M\to\infty}\sum_{k=1}^K\sum_{m=1}^M
(-1)^{m-1}(\frac\theta{m+\theta})^k.
\end{aligned}
\end{equation}
We have
\begin{equation}\label{KMlevel}
\begin{aligned}
&\sum_{k=1}^K\sum_{m=1}^M(-1)^{m-1}(\frac\theta{m+\theta})^k=\sum_{m=1}^M(-1)^{m-1}
\frac{\frac{\theta}{m+\theta}-(\frac{\theta}{m+\theta})^{K+1}}{1-\frac{\theta}{m+\theta}}=\\
&\theta\sum_{m=1}^M\frac{(-1)^{m-1}}m-\sum_{m=1}^M(-1)^{m-1}\thinspace\frac{m+\theta}m(\frac\theta{m+\theta})^{K+1}.
\end{aligned}
\end{equation}
Note that
\begin{equation}\label{log2}
\sum_{m=1}^\infty\frac{(-1)^{m-1}}m=\log 2.
\end{equation}
Since     $\frac{m+\theta}m(\frac\theta{m+\theta})^{K+1}$ is decreasing in $m$,
the second alternating series in \eqref{KMlevel}
satisfies the estimate
\begin{equation}\label{altser}
0\le \sum_{m=1}^M(-1)^{m-1}\thinspace\frac{m+\theta}m(\frac\theta{m+\theta})^{K+1}\le
(1+\theta)(\frac\theta{1+\theta})^{K+1}, \ \text{for}\ M,K\ge1.
\end{equation}
Now \eqref{toshowindep} follows from
\eqref{ontheway}-\eqref{altser}.\hfill $\square$

\end{document}